\documentclass[12pt]{article}
\usepackage[utf8]{inputenc}

\usepackage[top=1in, bottom=1in, left=1in, right=1in]{geometry}
\usepackage{graphicx} 
\usepackage{subfig}
\usepackage{amsthm}
\usepackage[authoryear]{natbib}
\usepackage{ulem}
\usepackage{color}
\usepackage{lineno}  
\usepackage{hyperref}
\usepackage{amsmath}
\usepackage{amssymb}
\usepackage{indentfirst} 
\usepackage{mathrsfs}

\numberwithin{equation}{section}
\newcommand{\eps}{\varepsilon}

\newcommand\keywords[1]{\textbf{Keywords}: #1}

\title{Bifurcations of a nonlinear spherical pendulum with vibrating suspension point}
\author{\small{Yan Luo$^2$ $\&$ Kaicheng Sheng$^{1,3}$\footnote{Corresponding author. E-mail addresses:  {\it k.sheng@sdu.edu.cn} (K.Sheng)}}\,, \\
\small{$^1$ Research Centre for Mathematics and Interdisciplinary Sciences, }\\
\small{Shandong University, Qingdao, 266237, China }\\
\small{$^2$ School of Mathematics and Statistics, Shandong University, Weihai, 264209, China} \\
\small{$^3$ Loughborough University, Loughborough, LE11 3TU, UK}
}

\date{}

\begin{document}

%\linenumbers

\maketitle

\begin{abstract}

This paper considers a nonlinear spherical pendulum whose suspension point performs high-frequency spatial vibrations. The dynamics of this pendulum can be described by averaging its Hamiltonian over phases of vibrations. Rotationally symmetric conditions on vibrations are assumed in the averaged Hamiltonian. Under these conditions, a bifurcation diagram for the phase portraits of the averaged system is presented. Numerical simulations of different examples of vibrations are performed. The case of proper degeneration in KAM theory guarantees the coherence of dynamical characteristics between the averaged and exact systems. %Bifurcation of phase portraits of the spherical physical pendulum with a vibrating suspension point is considered as well.

\bigskip

\noindent\keywords{nonlinear spherical pendulum, vibration, averaging method, phase portrait, KAM theory}
\end{abstract}

%\section*
%\,

%{\textbf{Small high frequency excitation has a considerable effect of dynamics of a mechanical system. A spherical pendulum is a classical model problem in mechanics along with a simple pendulum. Dynamics of a spherical pendulum with high frequency vertical harmonic vibration of the suspension point was considered by Markeyev (1999). We consider a spherical pendulum whose suspension point performs high-frequency arbitrary spatial vibrations, which has not been studied before. We construct averaged Hamiltonian and impose conditions on vibrations such that this Hamiltonian has a rotationally symmetry. We present a complete description of bifurcation diagram of its phase portraits of the averaged system. Bifurcations of phase portraits of spherical physical pendulum (rigid rod case) are discussed as well. We firstly give the full description of bifurcation diagram and its phase portraits of the averaged system.}}

\section{Introduction}

Small high-frequency excitation has a considerable effect on the dynamics of a mechanical system. \citet{steph_1, bogol_1,kapitsa_1} demonstrated it in a simple pendulum with a vertical vibrating suspension point. \citet{bogol_1} developed nonlinear theory by using the averaging method, and \citet{kapitsa_1} developed a method of separation of slow and fast variables for this (see also \citep{LL}). \citet{levi, bam, aio} and references therein considered the dynamics of a simple pendulum with a vibrating suspension point. It is demonstrated in \citep{AKN} that the problem can be simplified by averaging in Hamiltonian. The case of arbitrary planar vibrations of the suspension point of a planer simple pendulum is considered by \citet{nei}, which used the Hamiltonian approach of \citet{burd} to construct the averaged system and gave a complete description of bifurcations of its phase portraits. Generalisations to double-link and multiple-link pendulums are contained in \citep{steph_1, aches, khol}. 

The spherical pendulum, also known as the rigid body pendulum, is a classic example of a nonlinear dynamical system that exhibits rich and complex behaviour \citep{LL}. When the suspension point of the pendulum is fixed, its motion is constrained to a spherical surface and is governed by a set of coupled nonlinear differential equations. However, when the suspension point is subjected to vibrations, the dynamics of the pendulum become even more complicated, leading to a variety of nonlinear phenomena such as bifurcations, chaos, and resonance.

The study of bifurcations in the phase portraits of spherical pendulums with vibrating suspension points has garnered significant attention from researchers in the fields of nonlinear dynamics, control theory, and mechanical engineering. Bifurcations refer to qualitative changes in the behaviour of a dynamical system as a result of small changes in its parameters or initial conditions. In the context of spherical pendulums, it is interesting to consider bifurcations related to the emergence of new stable or unstable periodic orbits.

Several studies have investigated the dynamics of spherical pendulums with vibrating suspension points. For instance, \citet{mar} analyzed the stability and bifurcations of a spherical pendulum with a vertically oscillating suspension point analytically and showed that, by averaging over the fast vibration, the averaged system has either one or three equilibria in dependence on the system's parameters. \citet{malhotra1995stability} performed several numerical simulations and identified various bifurcation scenarios, including period-doubling bifurcations and the onset of chaos. Similarly, \citet{chen2007dynamics} studied the dynamics of a spherical pendulum with a horizontally oscillating suspension point and observed similar bifurcation phenomena. 

More recently, researchers have also explored the use of control strategies to stabilize or manipulate the dynamics of spherical pendulums with vibrating suspension points. For example, \citet{ortega2002stabilization} proposed a nonlinear control law to stabilize the pendulum's motion around an unstable equilibrium point, while \citet{chen2006chaotic} investigates the use of parametric excitation to induce chaotic dynamics in the system.

In this paper, a nonlinear spherical pendulum whose suspension point performs high-frequency vibrations in arbitrary directions is considered. A simple parametric form of the bifurcation curve in the plane of parameters is obtained. We construct the averaged Hamiltonian and impose conditions on vibrations such that the Hamiltonian has a rotationally symmetry. The relation between the exact and the averaged system is demonstrated by KAM theory. We aim to provide a comprehensive overview of the dynamics and bifurcations of spherical pendulums with vibrating suspension points. We present a complete description of the bifurcations of the phase portraits of the averaged system. Numerical examples demonstrate good agreement on the dynamics of the exact and the averaged system.

\section{Statement of the problem}
\label{Sec_Hamiltonian}

\begin{figure}[htbp] 
\centering
\includegraphics[width=0.3 \textwidth]{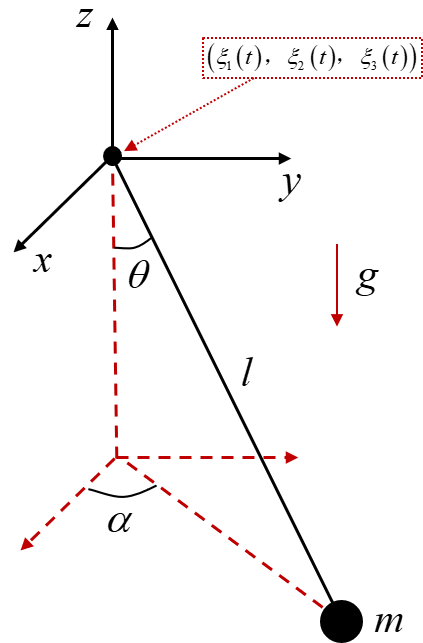} 
\caption{Nonlinear spherical pendulum}
\label{sphericalmodel} 
\end{figure} 

Considering a nonlinear spherical pendulum in $\mathbb{R}^3$, Fig.\ref{sphericalmodel}, whose suspension point performs high-frequency periodic vibrations.  $\xi_{i}\left(t \right)$ ($t\ge 0,\, i=1,2,3$) are Cartesian coordinates of the suspension point in the directions of $x$-axis, $y$-axis and $z$-axis, which demonstrate arbitrary vibrations of the suspension point. It is assumed that $\xi_{i}\left(t \right)$ are given periodic fast oscillating functions of time. Let $l=1$ and $m=1$ be the length of the massless rod and the mass of the bob for the pendulum. Denote by $\theta \in [0,\pi]$ the angle between the pendulum rod and the vertical line straight down ($-z$-axis) and by $\alpha \in [0, 2\pi]$ the angle of the rotation of the spherical pendulum. Owing to the spherical geometry of the problem, spherical coordinates are used to describe the position of the pendulum's bob:
\begin{equation}
\begin{aligned}
 & x=\sin \theta \cos \alpha +\xi_{1}\left(t\right),  \\ 
 & y=\sin \theta \sin\alpha +\xi_{2}\left(t\right),  \\ 
 & z=-\cos \theta +\xi_{3}\left(t\right). \\ 
\end{aligned} 
\end{equation}
The corresponding velocities are 
\begin{equation}
\begin{aligned}
 & \dot{x}=\cos \theta \cos \alpha \cdot \dot{\theta}-\sin \theta \sin \alpha \cdot \dot{\alpha} +\dot{\xi}_{1}\left( t \right),  
 \\ 
 & \dot{y}=\cos \theta \sin \alpha \cdot \dot{\theta}+\sin \theta \cos \alpha \cdot \dot{\alpha} +\dot{\xi}_{2}\left( t \right),  
 \\ 
 & \dot{z}=\sin \theta \cdot \dot{\theta }+\dot{\xi}_{3}\left( t \right). 
\end{aligned} 
\end{equation}
Then the kinetic energy of the bob up to some constants is 
\begin{equation}
\label{kinetic}
\begin{aligned}
T=&\frac{1}{2}m\left( {\dot{x}}^{2}+{\dot{y}}^{2}+{\dot{z}}^{2} \right)\\
=&\frac{1}{2}m{{l}^{2}}\left( {{{\dot{\theta }}}^{2}}
+{{\sin }^{2}}\theta \cdot {{{\dot{\alpha }}}^{2}} \right)+ml\left( \cos \theta \cos \alpha \cdot \dot{\theta }-\sin \theta \sin \alpha \cdot \dot{\alpha}  \right)\cdot \dot{\xi }\left( t \right) \\ 
&+ml\left( \cos \theta \sin \alpha \cdot \dot{\theta }+\sin \theta \cos \alpha \cdot \dot{\alpha } \right)\cdot \dot{\eta }\left( t \right)+ml\,\sin \theta \cdot \dot{\theta }\cdot \dot{\tau}\left( t \right), 
\end{aligned}
\end{equation}
and the potential energy with gravitational constant $g=1$ is
\begin{equation}
U=-\cos \theta +\xi_{3}(t).
\end{equation}
The generalised momenta conjugate to $\theta$ and $\alpha$ are
\begin{equation}
{{p}_{\theta }}=\frac{\partial {L}}{\partial \dot{\theta }}=\dot{\theta }+ \cos \theta \cos \alpha \cdot \dot{\xi}_{1}\left( t \right)+\cos \theta \sin \alpha \cdot \dot{\xi}_{2}\left( t \right)+\sin \theta \cdot \dot{\xi}_{3}\left( t \right) , 
\end{equation}
and 
\begin{equation}
{{p}_{\alpha }}=\frac{\partial {L}}{\partial \dot{\alpha }}={{\sin }^{2}}\theta \cdot \dot{\alpha } -\sin \theta \sin \alpha \cdot \dot{\xi}_{1}\left( t \right)+\sin \theta \cos \alpha \cdot \dot{\xi}_{2}\left( t \right),
\end{equation}
where ${L}=T-U$ is the Lagrangian of the system. Then the generalized velocities are
\begin{equation}
\label{velocities}
\begin{aligned}
  & \dot{\theta }={{p}_{\theta }}-\cos \theta \cos \alpha \cdot \dot{\xi}_{1}\left( t \right)-\cos \theta \sin \alpha \cdot \dot{\xi}_{2}\left( t \right)-\sin \theta \cdot \dot{\xi}_{3}\left( t \right),
  \\ 
 & \dot{\alpha }=\frac{{p}_{\alpha }}{{{\sin }^{2}}\theta}+\frac{\sin \alpha}{\sin \theta}\dot{\xi}_{1}\left( t \right)-\frac{\cos \alpha}{\sin \theta}\dot{\xi}_{2}\left( t \right). 
\end{aligned}
\end{equation}
The Hamiltonian for the spherical pendulum is a Legendre transformation of $L$: 
\begin{equation}
\label{original_H0}
H={p}_{\theta }\,\dot{\theta }+ {{p}_{\alpha }}\, \dot{\alpha }-{L}={p}_{\theta }\,\dot{\theta }+ {{p}_{\alpha }}\, \dot{\alpha }-T+U.
\end{equation}
Since $\theta$ and $\alpha$ are separated in $T$, we have 
\begin{equation}
T={{T}_{\theta, 1}}+{{T}_{\theta, 2}}+{{T}_{\alpha, 1}}+{{T}_{\alpha, 2}}, 
\end{equation}
where ${T}_{\theta, 1}$, ${T}_{\theta, 2}$ are terms of $\dot{\theta}$ and ${\dot{\theta }}^{2}$, and ${T}_{\alpha, 1}$, ${T}_{\alpha, 2}$ are terms of $\dot{\alpha}$ and ${\dot{\alpha}}^{2}$, respectively. Then 
\begin{equation}
\label{pthetapalpha}
\begin{aligned}
{p}_{\theta }\,\dot{\theta }&=\frac{\partial L}{\partial \dot{\theta }}\,\dot{\theta }=\frac{\partial T}{\partial \dot{\theta }}\,\dot{\theta }=\frac{\partial {{T}_{\theta, 2}}}{\partial \dot{\theta }}\,\dot{\theta }+\frac{\partial {{T}_{\theta, 1}}}{\partial \dot{\theta}}\,\dot{\theta }=2{T}_{\theta, 2}+{T}_{\theta, 1},
\\
{p}_{\alpha }\,\dot{\alpha }&=\frac{\partial L}{\partial \dot{\alpha }}\,\dot{\alpha }=\frac{\partial T}{\partial \dot{\alpha }}\,\dot{\alpha }=\frac{\partial {{T}_{\alpha, 2}}}{\partial \dot{\alpha }}\,\dot{\alpha}+\frac{\partial {{T}_{\alpha, 1}}}{\partial \dot{\alpha }}\,\dot{\alpha }=2{T}_{\alpha, 2}+{T}_{\alpha, 1}.
\end{aligned}
\end{equation}
Therefore, from formula (\ref{original_H0}) and (\ref{pthetapalpha}), we have
\begin{equation}
\begin{aligned}
H&=2\,{T}_{\theta, 2}+{T}_{\theta, 1}+2\,{T}_{\alpha, 2}+{T}_{\alpha, 1}-\left({{T}_{\theta, 1}}+{{T}_{\theta, 2}}+{{T}_{\alpha, 1}}+{{T}_{\alpha, 2}}\right)+U
\\
&={T}_{\theta, 2}+{T}_{\alpha, 2}+U=\frac{1}{2}m{{l}^{2}}\left( {{{\dot{\theta }}}^{2}}
+{{\sin }^{2}}\theta \cdot {{{\dot{\alpha }}}^{2}} \right)-\cos \theta +\xi_{3}\left( t \right).
\end{aligned}
\end{equation}
By formula (\ref{velocities}), the Hamiltonian turns to be
\begin{equation}
\label{original_H}
\begin{aligned}
H=&\frac{1}{2}{{\left[{{p}_{\theta }}-\cos \theta \cos \alpha \cdot \dot{\xi}_{1}\left( t \right)-\cos \theta \sin \alpha \cdot \dot{\xi}_{2}\left( t \right)-\sin \theta \cdot \dot{\xi}_{3}\left( t \right) \right]}^{2}} \\ 
&+\frac{1}{2\,{{{\sin }^{2}}\theta } }{ \left[{{p}_{\alpha }}+\sin \theta \sin \alpha \cdot \dot{\xi}_{1}\left( t \right)-\sin \theta \cos \alpha \cdot \dot{\xi}_{2}\left( t \right)\right]^{2} }-\cos \theta +\xi_{3}\left( t \right). 
\end{aligned}
\end{equation}

\section{Averaged Hamiltonian of the system}

Assume that 
\begin{equation}
\xi_{i}\left(t \right)=\varepsilon \tilde{\xi}_{i}\left(\omega t/\varepsilon\right), \quad i=1,2,3.
\end{equation}
where $ \varepsilon $ is a small parameter, $\tilde{\xi}_{i}$ are $2\pi\varepsilon/\omega$-periodic functions of $t$ with zero average value. Then $\dot{\xi}_{i}$ are values of order $1$. In the system of Hamilton's equations:  
\begin{equation}
\dot{\theta} = \frac{\partial H}{\partial p_{\theta}}, \quad \dot{p_{\theta}} = - \frac{\partial H}{\partial \theta}; \quad  \dot{\alpha} = \frac{\partial H}{\partial p_{\alpha}}, \quad \dot{p_{\alpha}} = - \frac{\partial H}{\partial \alpha}.
\end{equation}
the right-hand sides are fast oscillating functions of time. {The dynamics} of variables $\theta$, $p_{\theta}$, $\alpha$, $p_{\alpha}$ can be approximately described by averaging the Hamiltonian $H$ with respect to time $t$ over {the phases} of vibrations \citep{BM}. The averaged value of the vibrations are:
\begin{equation}
\begin{aligned}
\overline{{{\xi}_{i}}\left( t \right)}&=\frac{1}{2\pi}\int_{0}^{2\pi}{{\xi}_{i}\left( t \right)}\,d{\,\frac{\omega t}{\varepsilon}} =0;
\\
\overline{{\dot{\xi}_{i}}\left( t \right)}&=\frac{1}{2\pi}\int_{0}^{2\pi}{{\dot{\xi}_{i}\left( t \right)}}\,d{\,\frac{\omega t}{\varepsilon}} =0;
\\
\overline{{\dot{\xi}_{i}}\left( t \right){\dot{\xi}_{j}}\left( t \right)}&=\frac{1}{2\pi}\int_{0}^{2\pi}{{\dot{\xi}_{i}}\left( t \right){\dot{\xi}_{j}}\left( t \right)}\,d{\,\frac{\omega t}{\varepsilon}}, \quad  i,\,j=1,2,3.
\end{aligned}
\end{equation}
Then the averaged Hamiltonian of the spherical pendulum system is
 %Consider the system of equations given in action-angle variables $\mathbf{I}=(I_1, \ldots, I_n)$, $\boldsymbol{\theta}= (\theta_1, \ldots, \theta_n)\  {\rm mod\,2\pi}$
%\begin{equation*}
%\begin{aligned}
%\dot{\mathbf{I}}&=\eps\,\mathbf{g}\left(\mathbf{I},\boldsymbol{\theta}\right),
%\\
%\dot{\boldsymbol{\theta}}&=\boldsymbol{\omega}\left(\mathbf{I}\right)+\eps\,\mathbf{f}\left(\mathbf{I},\boldsymbol{\theta}\right).
%\end{aligned}
%\end{equation*}
%$\mathbf{g}$ and $\mathbf{f}$ are $2\pi$-periodic functions of $\boldsymbol{\theta}$. The corresponding averaged system
%\begin{equation*}
%\dot{\mathbf{J}}=\eps\,\bar{\mathbf{g}}\left(\mathbf{J}\right),\, \bar{\mathbf{g}}\left(\mathbf{J}\right)={\left({1}/{2 \pi}\right)}^{n}\,\int\nolimits_{0}^{2\pi}\ldots\int\nolimits_{0}^{2\pi}{\mathbf{g}\left(\mathbf{J},\theta\right)}d{\theta_{1}} \ldots d{\theta_{n}}.
%\end{equation*}
%The difference between the slow variables $\mathbf{I}$ in exact system and $\mathbf{J}$ in averaged system is small over the times of order $1/\eps$: 
%\begin{equation*}
%\left|\mathbf{I}\left(t\right)-\mathbf{J}\left(t\right)\right|<c\,\eps,\,\,\,if \,\,\,\mathbf{I}\left(0\right)=\mathbf{J}\left(0\right), \,\,\,0\le t \le 1/\eps, \,\,c=\text{const}.
%\end{equation*}
%($2\pi$-periodic functions of the argument $\omega_{i} t/\varepsilon$). 
\begin{equation}
\label{averagedH}
\begin{aligned}
\bar{H}&=\frac{1}{2}\left( {{p}_{\theta }}^{2}+\frac{{{p}_{\alpha }}^{2}}{{{\sin }^{2}}\theta } \right)+\frac{1}{2}\left( {{\cos }^{2}}\theta \, {{\cos }^{2}}\alpha +{{\sin }^{2}}\alpha  \right)\cdot \overline{{\dot{\xi}_{1}^{2}}\left( t \right)} \\ 
 & +\frac{1}{2}\left( {{\cos }^{2}} \theta\,  {{\sin }^{2}}\alpha +{{\cos }^{2}}\alpha  \right)\cdot \overline{{\dot{\xi}_{2}^{2}}\left( t \right)} +\frac{1}{2}{{\sin }^{2}}\theta \cdot \overline{{\dot{\xi}_{3}^{2}}\left( t \right)} \\ 
 & +\left( {{\cos }^{2}} \theta \cos \alpha \sin \alpha -\cos \alpha \sin \alpha  \right)\cdot \overline{\dot{\xi}_{1}\left( t \right)\dot{\xi}_{2}\left( t \right)}+\cos \theta \sin \alpha \sin \theta \cdot \overline{\dot{\xi}_{2}\left( t \right)\dot{\xi}_{3}\left( t \right)}\\ 
 & +\cos \theta \cos \alpha \sin \theta \cdot \overline{\dot{\xi}_{3}\left( t \right)\dot{\xi}_{1}\left( t \right)} -\cos \theta .
\end{aligned}
\end{equation}
%Here ``bars'' on the right-hand side denote averaging over time.
%\section{The Hamiltonian under symmetry  conditions}

The averaged Hamiltonian is in the absence of $\alpha$ if vibrations $\xi_{i}\left(t \right)$ satisfy the rotationally symmetric conditions
\begin{equation}
\label{symmetry}
\begin{aligned}
  & a) \quad \overline{\dot{\xi_{1}^{2}}\left( t \right)}=\overline{\dot{\xi_{2}^{2}}\left( t \right)}, \\ 
 & b) \quad \overline{\dot{\xi_{i}}\left( t \right)\dot{\xi_{j}}\left( t \right)}=0, \quad i,j=1,2,3, \,i\neq j. 
\end{aligned}
\end{equation}
With these conditions, $\alpha$ is a cyclic coordinate, and $p_{\alpha}$ is a first integral of the averaged system. Then the averaged Hamiltonian (\ref{averagedH}) simplifies to
\begin{equation}
\label{av_ham0}
\bar{H}=\frac{{{p}_{\theta }}^{2}}{2}+\frac{{{p}_{\alpha}}^{2}}{2 {{\sin }^{2}}\theta }+\frac{1}{2}\, {{\cos }^{2}}\theta\cdot \overline{\dot{\xi_{2}^{2}}\left( t \right)}+\frac{1}{2}\,\overline{\dot{\xi_{2}^{2}}\left( t \right)}+\frac{1}{2}\,{{\sin }^{2}}\theta \cdot \overline{\dot{\xi_{3}^{2}}\left( t \right)}-\cos \theta.
\end{equation}
The term $\frac{1}{2}\,\overline{\dot{\xi_{2}^{2}}\left( t \right)}$ can be omitted because it is a constant that does not change the stability of the spherical pendulum. For convenience, we take $A=\overline{\dot{\xi_{3}^{2}}\left( t \right)}$, $B={{p}_{\alpha }}^{2}$ and $C=\overline{\dot{\xi_{1}^{2}}\left( t \right)}=\overline{\dot{\xi_{2}^{2}}\left( t \right)}$. Up to a constant, equation (\ref{av_ham0}) turns to
\begin{equation}
\label{av_ham}
\bar{H}=\frac{{{p}_{\theta }}^{2}}{2}+\bar{U},
\end{equation}
where the effective potential
\begin{equation}
\bar{U}=\frac{B}{2{{\sin }^{2}}\theta }+\frac{1}{2}{{\sin }^{2}}\theta \cdot (A-C)-\cos \theta.
\end{equation}
The dynamics of variables $\theta$, $p_{\theta}$ are described by the 1 degree of freedom (1-DOF) Hamiltonian system $\bar{H}$. By equation (\ref{velocities}), the behaviour of the variable $\alpha$ is described by the equation
\begin{equation}
\dot{\alpha }=\frac{{p}_{\alpha }}{{{\sin }^{2}}\theta}+\frac{\sin \alpha}{\sin \theta}\cdot\,\overline{\dot{\xi}_{1}\left( t \right)}-\frac{\cos \alpha}{\sin \theta}\cdot\,\overline{\dot{\xi}_{2}\left( t \right)}=\frac{p_{\alpha}}{{\sin }^2\theta }. 
\end{equation}

\begin{figure}[htbp] 
  \centering
\subfloat[]{\includegraphics[width=0.45\textwidth]{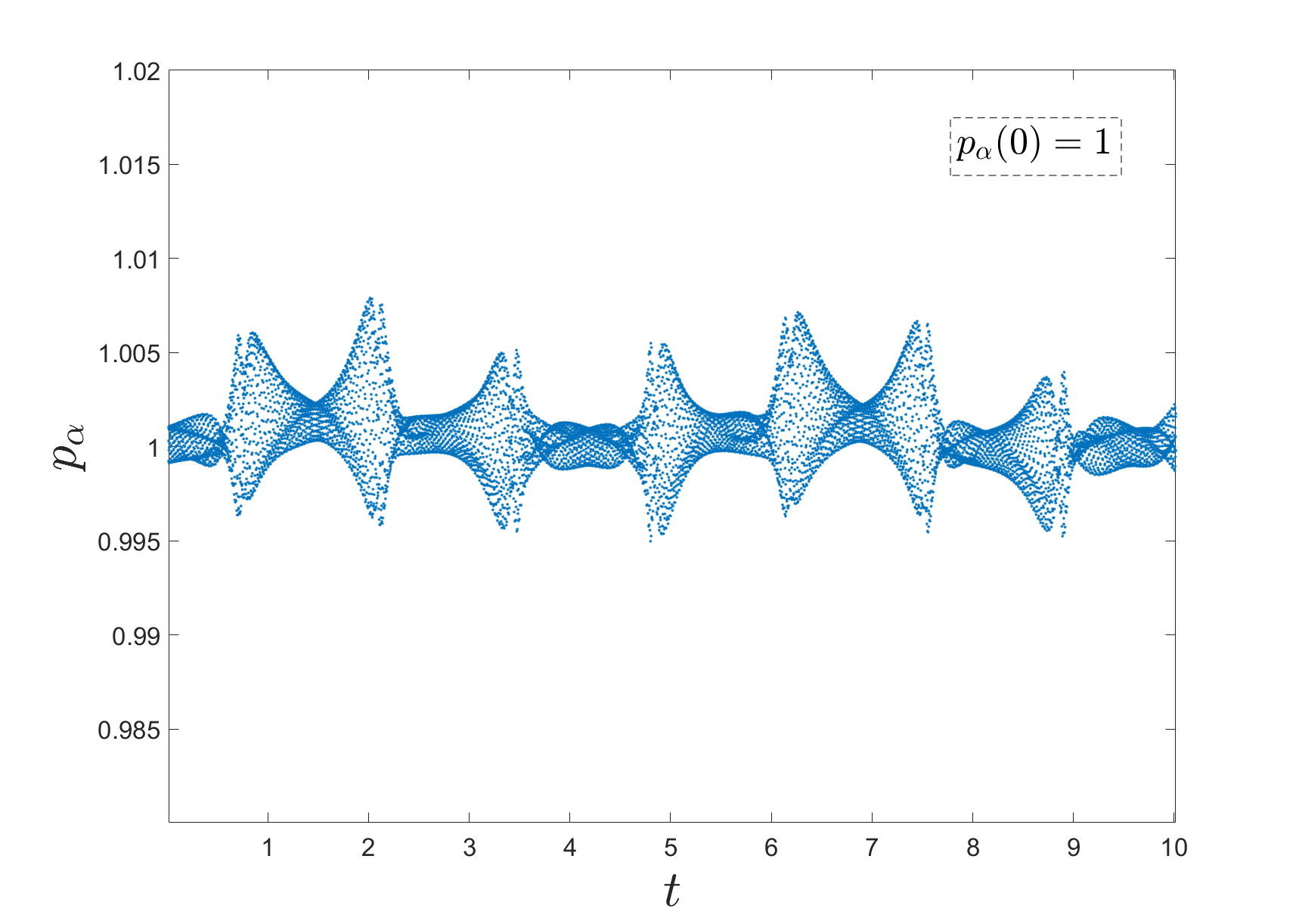}\label{graph_p_alpha1}} 
\subfloat[]{\includegraphics[width=0.45\textwidth]{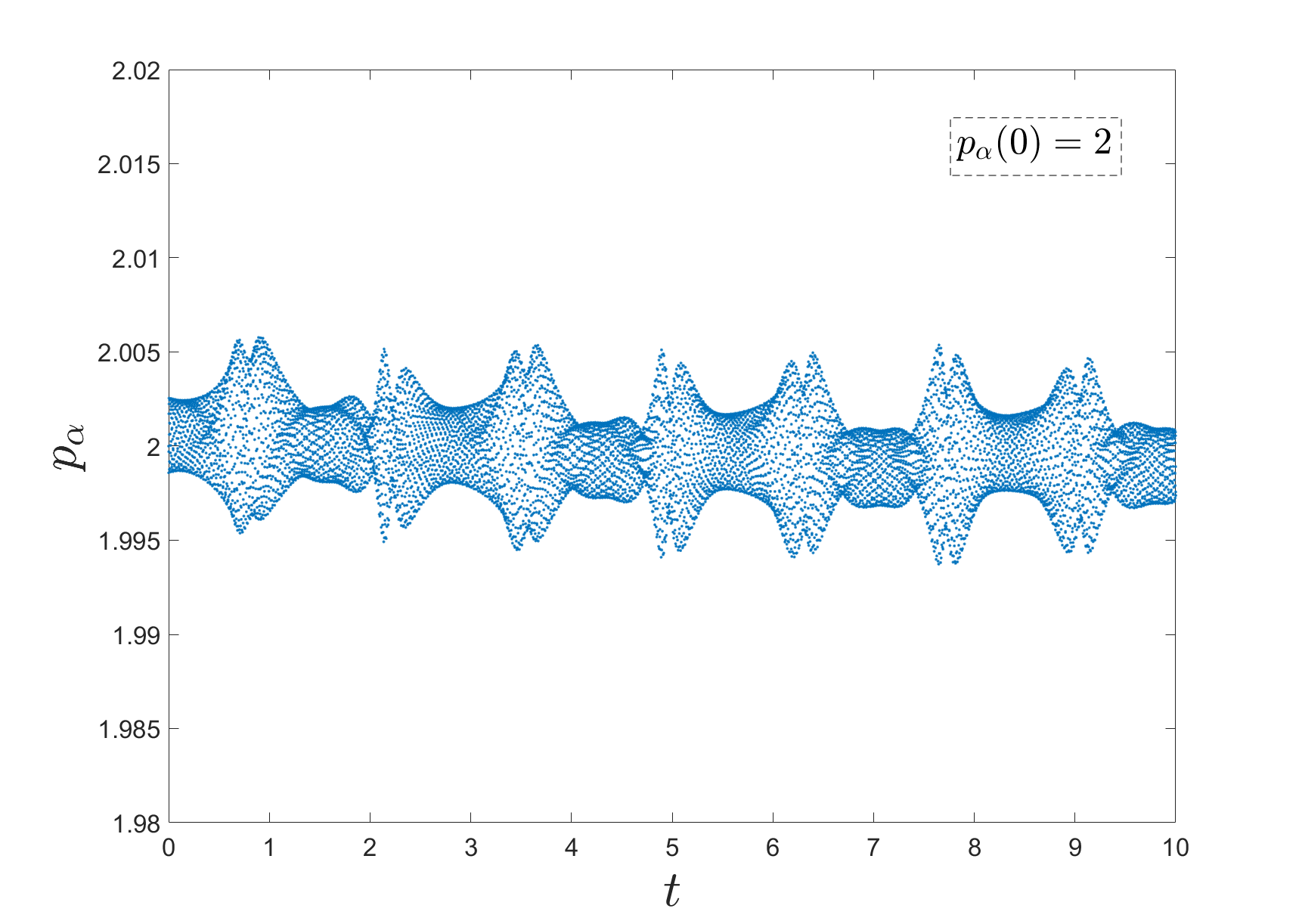}\label{graph_p_alpha3}} 
\caption{Behaviour of $p_{\alpha}$ over $t$ in exact system.}
%\label{graph_p_alpha}
\end{figure}

Taking $\xi_{1}=\varepsilon\sin \left( {t}/{\varepsilon}\right)$, $\xi_{2}=\varepsilon\cos \left( {t}/{\varepsilon}\right)$ and $\xi_{3}=3\varepsilon\sin \left( {2 t}/{\varepsilon}\right)$ as an example of the vibrations, one can verify that these vibrations satisfy the conditions (\ref{symmetry}). Numerical simulations of the behaviour of $p_{\alpha}$ over $t$ in the exact (not averaged) system are performed for these values with initial value $p_{\alpha}(0)=1$ (Fig. \ref{graph_p_alpha1}) and $p_{\alpha}(0)=2$ (Fig. \ref{graph_p_alpha3}). The values of $p_{\alpha}(t)$ move around the initial value, thus $p_{\alpha}(t)$ is close to a first integral. The results match the averaged results well.

\section{Relation to exact problem}

Under the rotationally symmetric conditions (\ref{symmetry}), in the 4-dimensional phase space of the averaged system $\bar{H}\left(p_{\theta}, \theta, p_{\alpha}, \alpha \right) $, the motion occurs on invariant surfaces $p_{\alpha}=c_1, \bar{H}=c_2$, $c_1$, $c_2$ are some constants. These surfaces could be 2-dimensional tori, closed trajectories, or separatrix surfaces. We consider here a compact part of the phase space and assume that vibrations $\xi_{i}\left( t \right)$ ($i=1,2,3$) are described by smooth enough functions ($\xi_{i}\in {\mathcal{C}}^k \left( \mathbb{R} \right)$, $k$ is large enough), then, KAM theory (see, e.g.,  \cite{AKN}) is applicable. 

By introducing a new time $\tilde{t}=t / \eps$ , we obtain a new Hamiltonian system of 3-DOF: 
\begin{equation}
\hat{H}={I}_{\tilde{t}}+H\left(p_{\theta}, \theta, p_{\alpha}, \alpha, \tilde{t} \right) 
\end{equation}
with $\dot{\tilde{t}}={\partial \hat{H}} / {\partial {I}_{\tilde{t}}}=1$. Then we consider the Hamiltonian
\begin{equation}
\tilde{H}=\eps \hat{H}=\tilde{I}_{\tilde{t}}+\eps H\left(p_{\theta}, \theta, p_{\alpha}, \alpha, \tilde{t} \right) ,
\end{equation}
in which the variables $\theta,\alpha, p_{\theta}, p_{\alpha}$ are slow and the variable $\tilde{t}$ is fast, $\tilde{I}_{\tilde{t}}$ and $\tilde{t}$ are conjugate variables. We generated a 2.5-DOF Hamiltonian system $H$ to a 3-DOF Hamiltonian system $\tilde{H}$. The term  $\tilde{I}_{\tilde{t}}$ can be regarded as a properly degenerate Hamiltonian because 
\begin{equation}
\det\left({{{\partial }^{2}}{{\tilde{I}_{\tilde{t}}}}} / {\partial {{\left(p_{\theta}, p_{\alpha}, \tilde{I}_{\tilde{t}} \right)}^{2}}}\right)\equiv 0.
\end{equation}

The Hamiltonian $\eps H$ (or the averaged Hamiltonian $\eps \bar{H}$) is non-degenerate with respect to $p_{\theta}$, $p_{\alpha}$, which means the perturbation removes the degeneracy of the system. Then according to KAM theory, the 5-dimensional extended phase space $\left(\theta,\alpha, p_{\theta}, p_{\alpha}, t \right)$ of the original system is filled by the 3-dimensional invariant tori close to the tori
$$
\bar{H}=c_1, \quad p_{\alpha}=c_2, \quad c_1, c_2={\rm const}.
$$
The remainder is of a measure $O(\exp (-{\rm const}/\eps)$, which is exponentially small when $\varepsilon$ is small \citep{AKN}. This is the case of proper degeneration in KAM theory. 

If rotationally symmetric conditions (\ref{symmetry}) are satisfied only approximately, with some accuracy $\delta$, then KAM theory ensures that the 5-dimensional extended phase space of the original system is filled by 3-dimensional invariant tori close to tori $p_{\alpha}=c_1, \bar{H}=c_2$ up to a reminder of a measure which is small when $\delta$ and $\varepsilon$ are small. 

Consequently, the dynamics of the exact Hamiltonian system of the spherical pendulum are close to the averaged Hamiltonian system under the rotationally symmetric conditions or close to them. We will analyse the dynamics of the averaged Hamiltonian system in the following.

\section{Bifurcations of phase portraits}

We will divide the parameter plane of the problem into two domains corresponding to different types of phase portraits of the averaged system. The boundary between these domains is a bifurcation curve corresponding to degenerate equilibria: first and second derivatives of $\bar{U}$ vanish for parameters on these curves. The number of equilibria changes at crossing such a curve in the plane of parameters.

The bifurcation curve is defined by the equations
\begin{equation}
\label{parametric_1}
\begin{aligned}
\frac{\partial \bar{U}}{\partial \theta}=&-\frac {B\cos \theta}{ {{\sin }^{3}}\theta}+ \frac{1}{2} (A-C)\sin 2 \theta +\sin\theta =0,
 \\
\frac{{{\partial }^{2}}\bar{U}}{\partial {{\theta}^{2}}}=&\frac {3 B\,{{\cos}^{2}}\theta}{ {{\sin }^{4}}\theta}+\frac {B}{{\sin }^2\theta }+ (A-C) \cos 2\theta +\cos\theta=0.
\end{aligned} 
\end{equation}

\begin{figure}[htbp!]
  \centering
\subfloat[]{\includegraphics[width=0.45\textwidth]{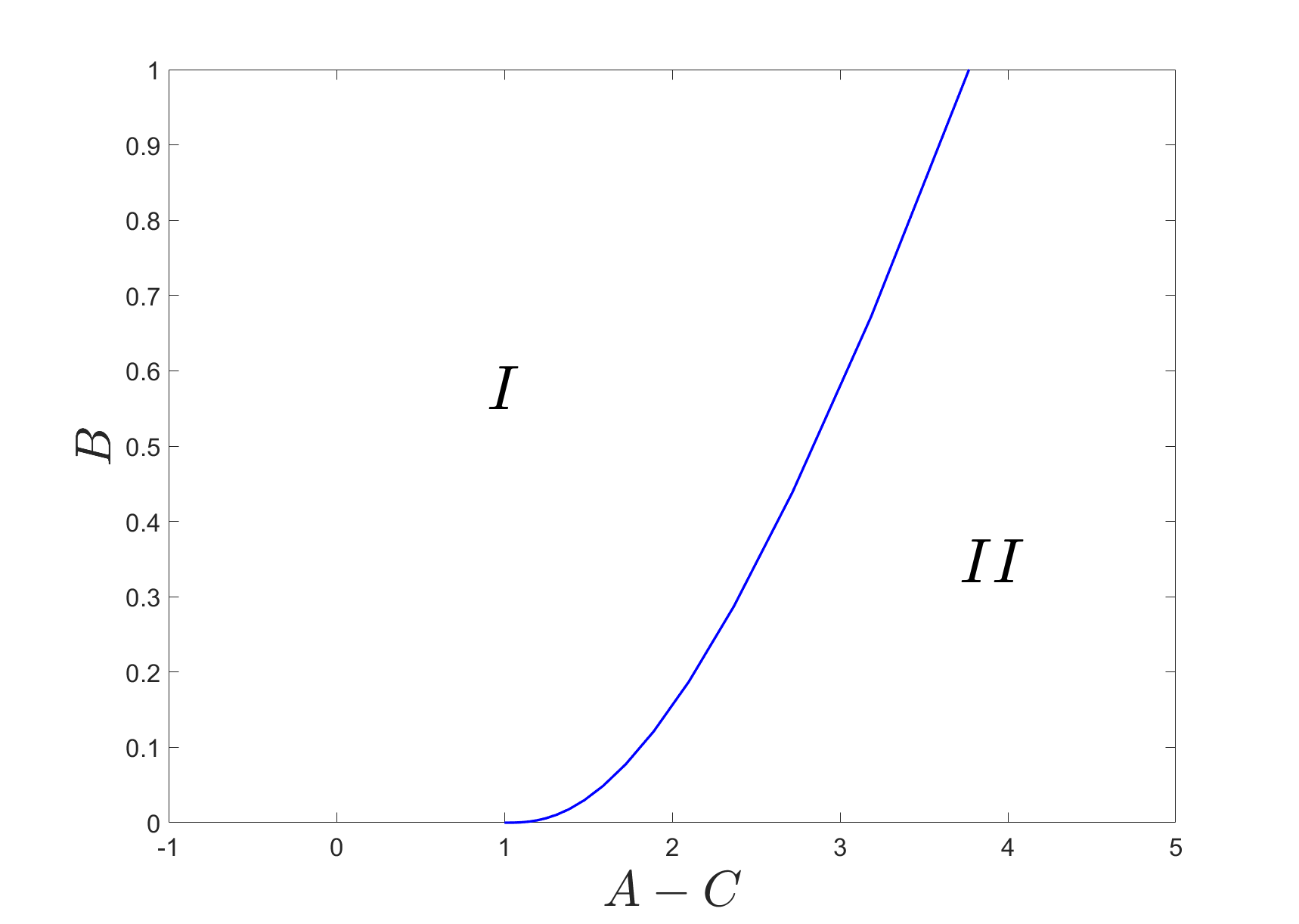}\label{critical1}} 
\subfloat[]{\includegraphics[width=0.45\textwidth]{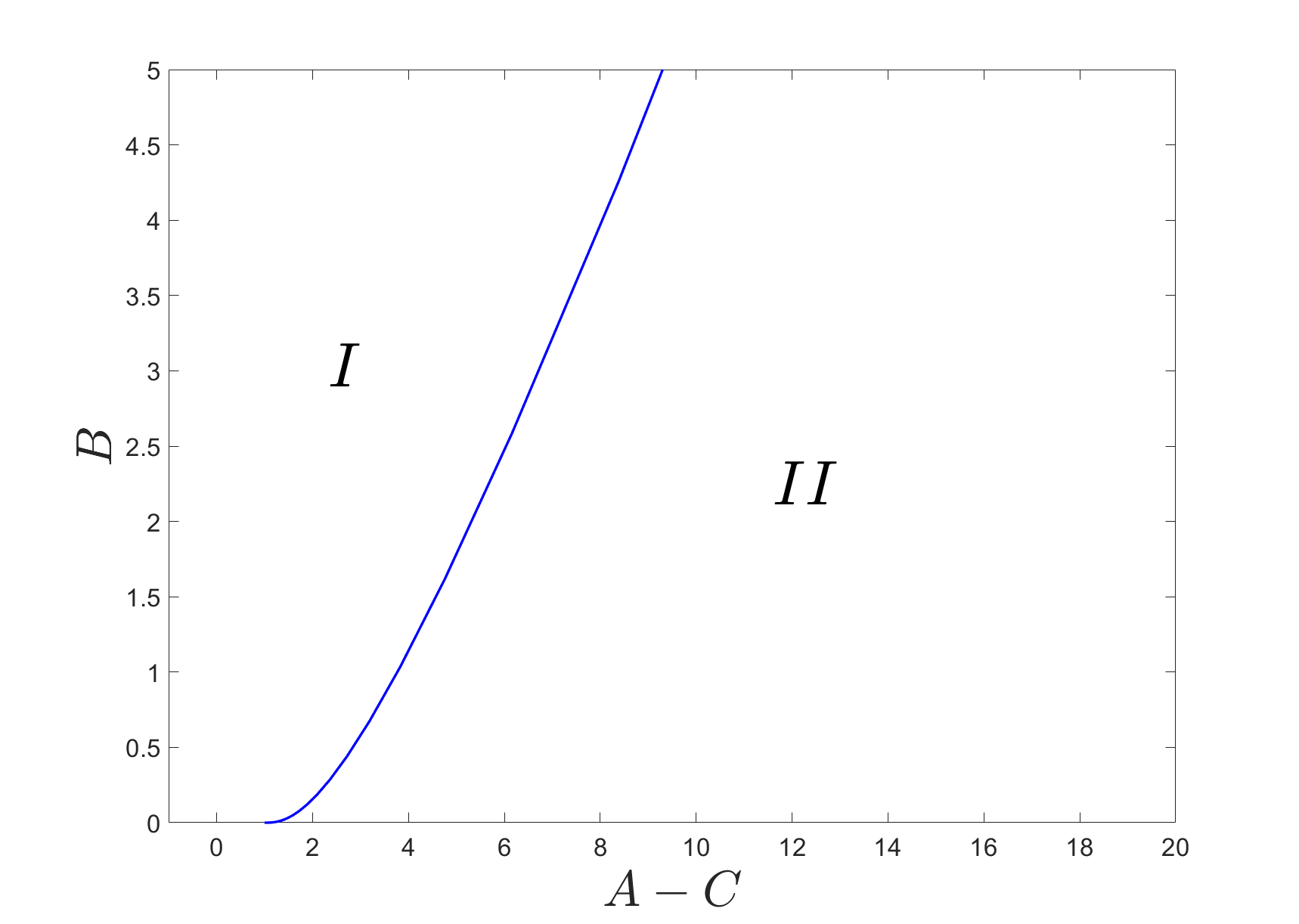}\label{critical2}} 
\caption{Bifurcation curve corresponds to degenerate equilibria.}
\label{fig:2} 
\end{figure}

\begin{figure}[htbp!]
  \centering
\subfloat[Domain $I$]{\includegraphics[width=0.45\textwidth]{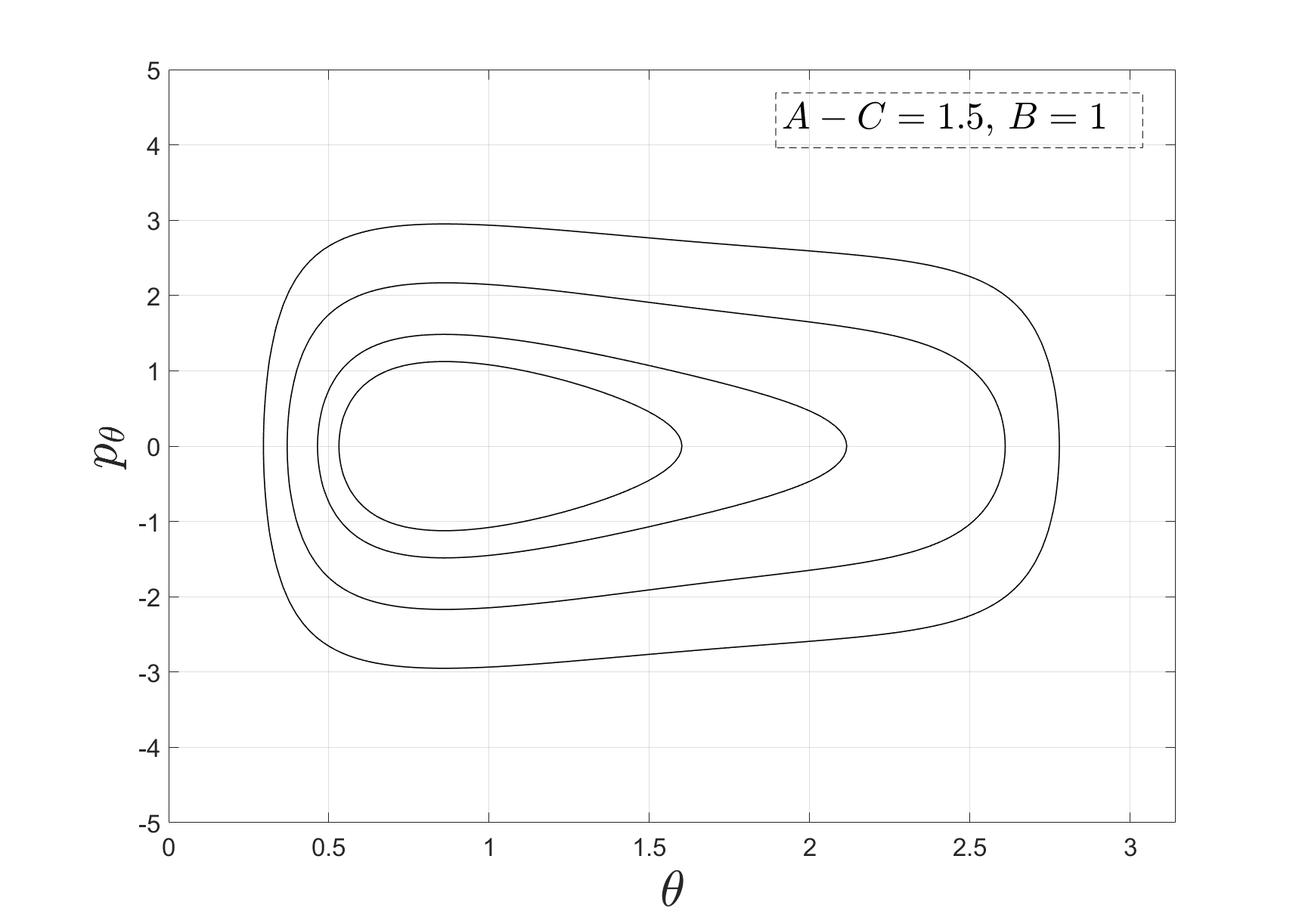}\label{phase3}} 
\subfloat[Domain $II$]{\includegraphics[width=0.45\textwidth]{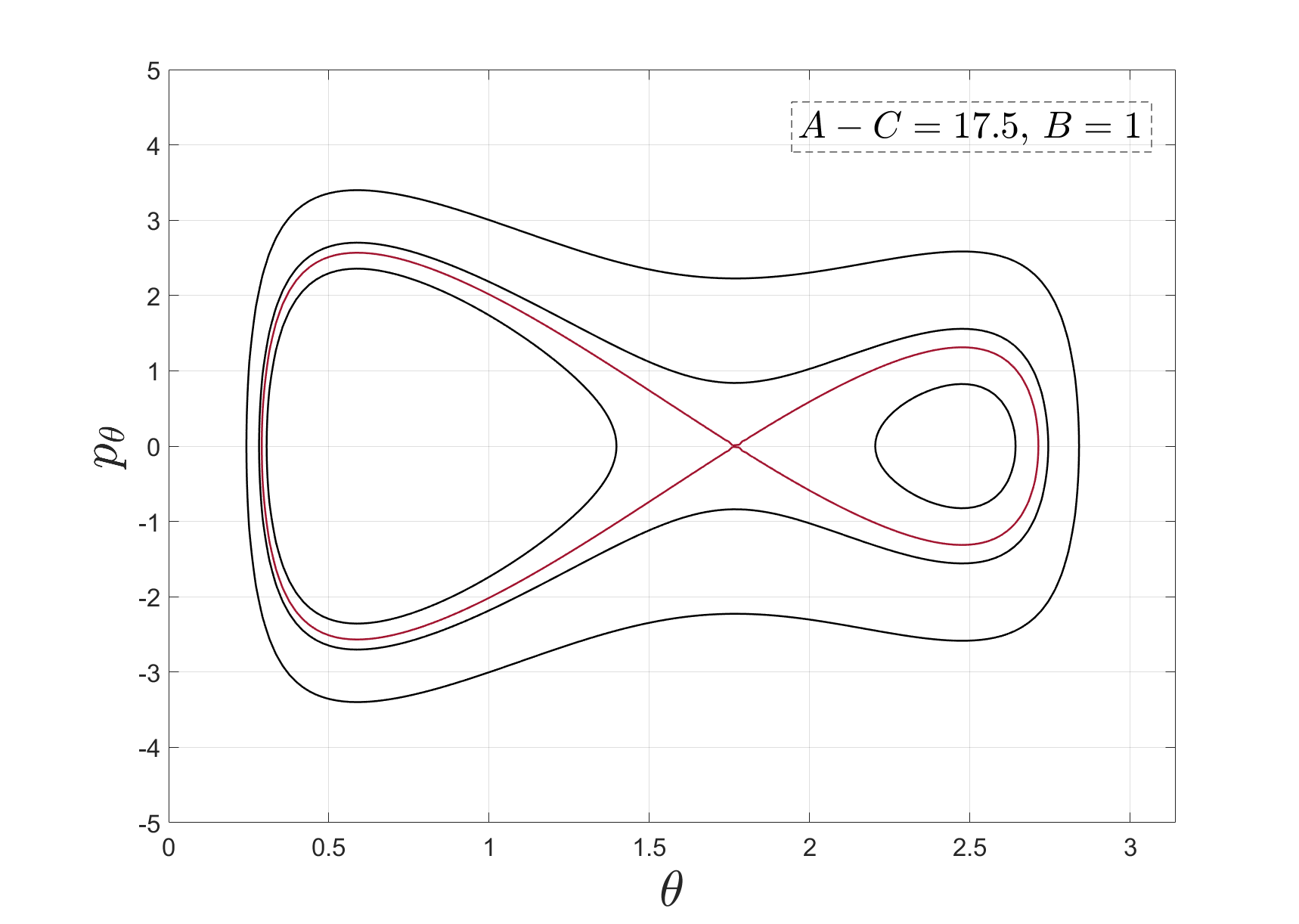}\label{phase4}} 
\caption{Phase portrait for parameters.}
%\label{fig:3} 
\end{figure}

Equations (\ref{parametric_1}) is a system of linear non-homogeneous equations for $A-C$ and $B$. Solving it, the parametric representation of the bifurcation curve with the parameter $\theta$ is
\begin{equation}
\label{parametric}
A-C=-\frac{3{{\cos}^{2}}\theta +1}{4{{\cos }^{3}}\theta }, \quad B=-\frac{1}{4}\cdot \frac{{{\sin }^{6}}\theta }{{{\cos }^{3}}\theta }. 
\end{equation}
As $B={{p}_{\alpha }}^{2}$ is non-negative, then $\theta \in [\pi/2,\pi]$. Therefore there are no bifurcations with $\theta \in [0,\pi/2)$. The bifurcation curve $\Gamma$ corresponding to degenerate equilibria is shown in Fig.\ref{critical1} and \ref{critical2} in different scales. It divides the plane $(A-C,\, B)$ into domains $I$ and $II$. The line $B=0$ corresponds to a simple pendulum (cf. \cite{nei}). Phase portraits of the system for parameters in domains $I$ and $II$ are shown in Fig.\ref{phase3} and \ref{phase4} respectively. The averaged spherical pendulum has one equilibrium if the parameters are in domain $I$, and three equilibria if the parameters are in domain $II$. 

\section{On the Poincar\'e return map}

In the plane of $\theta$, $p_{\theta}$, we consider the Poincar\'e return map with Poincar\'e section $t=0 \ {\rm mod}\, \frac{2\pi\eps}{\omega}$. There are fixed points, invariant curves and manifolds of hyperbolic fixed points in the Poincar\'e return map.  The fixed points of this map are close to the equilibria of the averaged system. Invariant curves filled the $\theta-p_{\theta}$ plane up to a remainder of a measure which is exponentially small in the perturbation parameter $\eps$.  Stable and unstable manifolds of the hyperbolic fixed point in the return map are split with exponentially small distances, the detailed estimation of the splitting is an open question. Stable and unstable manifolds in the return map are close to separatrices of the averaged system. 

Thus, our study of phase portraits of the averaged system is closely related to the dynamics of the exact (not averaged) problem. The dynamical characteristics of the exact problem are similar to the dynamical characteristics in the averaged system. 

\section{Numerical examples}

In this section, we provide different numerical examples with suspension point vibrations. 

\begin{table}[htbp!]
\caption{Examples of vibrations}
\label{table1}
    \centering
    \begin{tabular}{|l|l|l|l|l|l|l|}
    \hline
         & $\xi_{1}$ & $\xi_{2}$ & $\xi_{3}$ & $A-C$ & $p_{\alpha}(0)$ \\ \hline
         Example 1 & $\varepsilon\sin \left( {t}/{\varepsilon}\right)$ & $\varepsilon\cos \left( {t}/{\varepsilon}\right)$ & $3\varepsilon\sin \left( {2 t}/{\varepsilon}\right)$ & 17.5 & 1 \\ \hline
         Example 2 &  $\varepsilon\sin \left( { t}/{\varepsilon}\right)$ &  $\varepsilon\cos \left( { t}/{\varepsilon}\right)$ & $\varepsilon\sin \left( { 2 t}/{\varepsilon}\right)$ & 1.5 & 1 \\ \hline
         Example 3 & $\varepsilon\sin \left( {t}/{\varepsilon}\right)$ & $\varepsilon\cos\left( {t}/{\varepsilon}\right)$ & $3\varepsilon\sin \left( {2 t}/{\varepsilon}\right)$ & 17.5 & 2\\ \hline  
         Example 4 &  $\varepsilon\sin \left( { t}/{\varepsilon}\right)$ &  $\varepsilon\cos \left( { t}/{\varepsilon}\right)$ & $\varepsilon\sin \left( { 2 t}/{\varepsilon}\right)$ & 1.5 & 2 \\ \hline
    \end{tabular}
\end{table}

Rotationally symmetric conditions (\ref{symmetry}) are satisfied for these examples of vibrations. We take $\varepsilon=0.001$. Then Example 1 and Example 2 correspond to points $P(17.5,1)$ in domain II and  $Q(1.5,1)$ in domain I in the parameter plane respectively with initial value $p_{\alpha}(0)=1$, Example 3 and Example 4 correspond to points $P(17.5,4)$ in domain II and  $Q(1.5,4)$ in domain I in the parameter plane respectively with initial value $p_{\alpha}(0)=2$. 

\begin{figure}[htbp!]
  \centering
\subfloat[]{\includegraphics[width=0.45\textwidth]{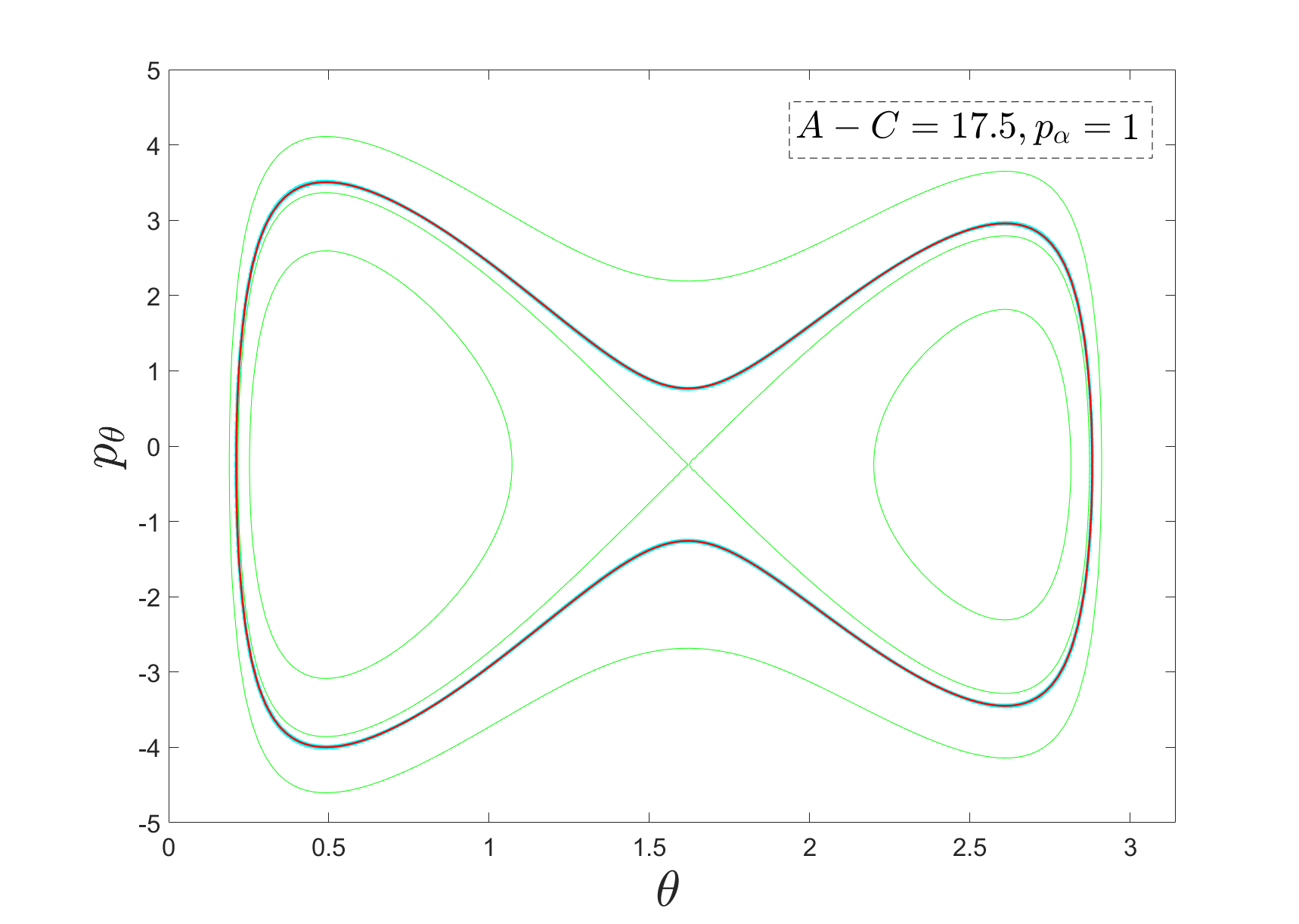}\label{p_phi_1}} 
\subfloat[]{\includegraphics[width=0.45\textwidth]{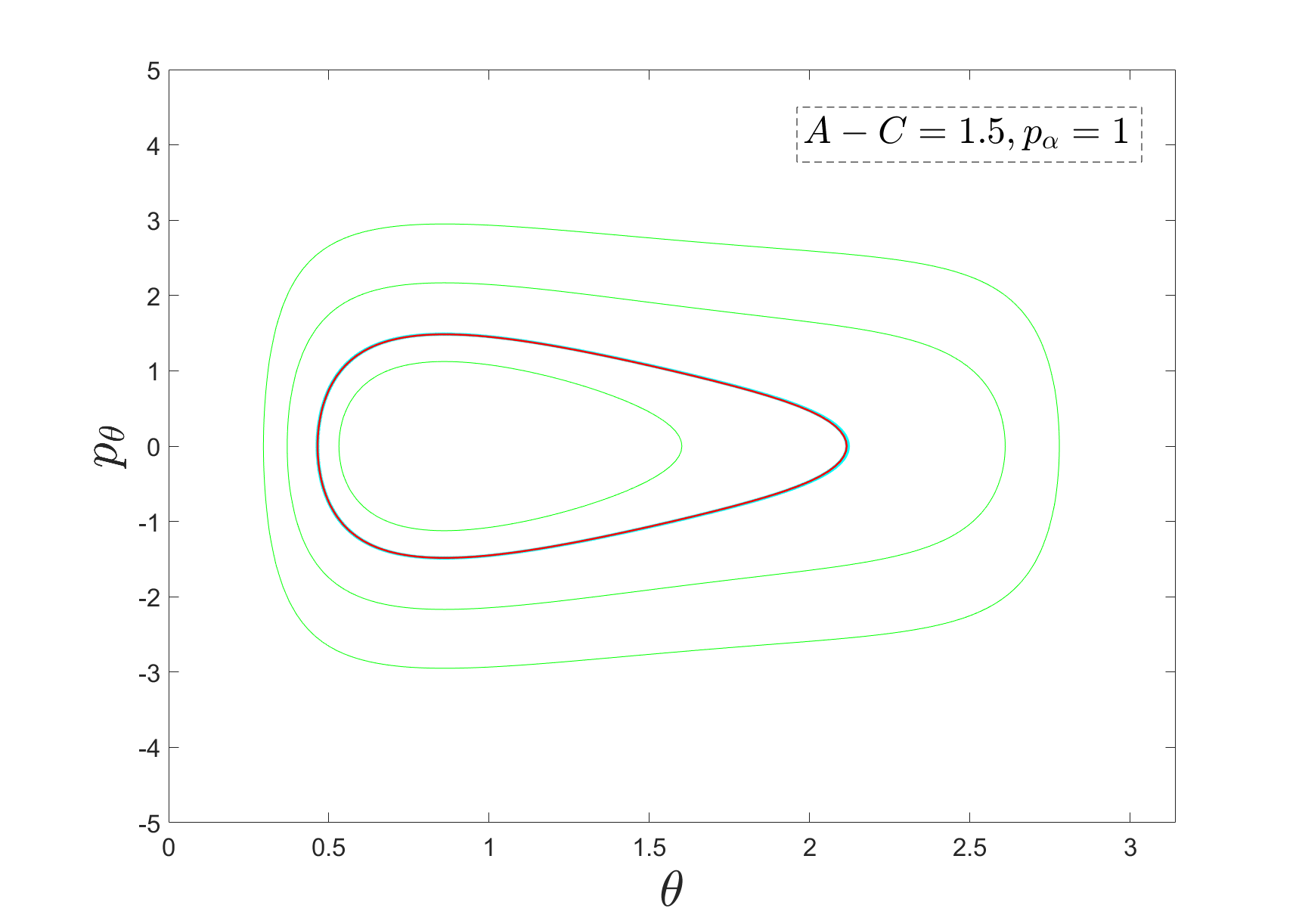}\label{p_phi_2}} 
\\
\subfloat[]{\includegraphics[width=0.45\textwidth]{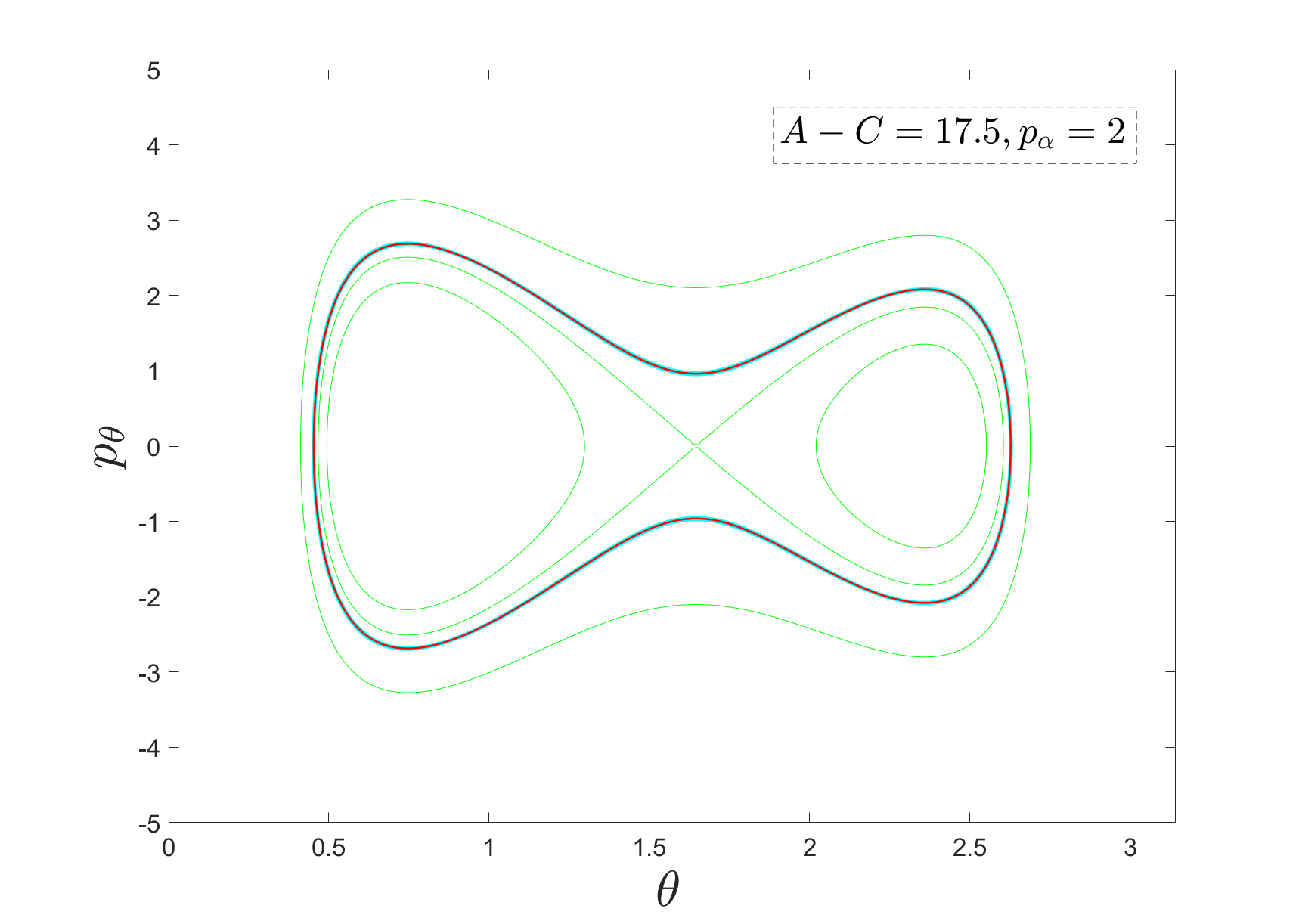}\label{p_phi_3}} 
\subfloat[]{\includegraphics[width=0.45\textwidth]{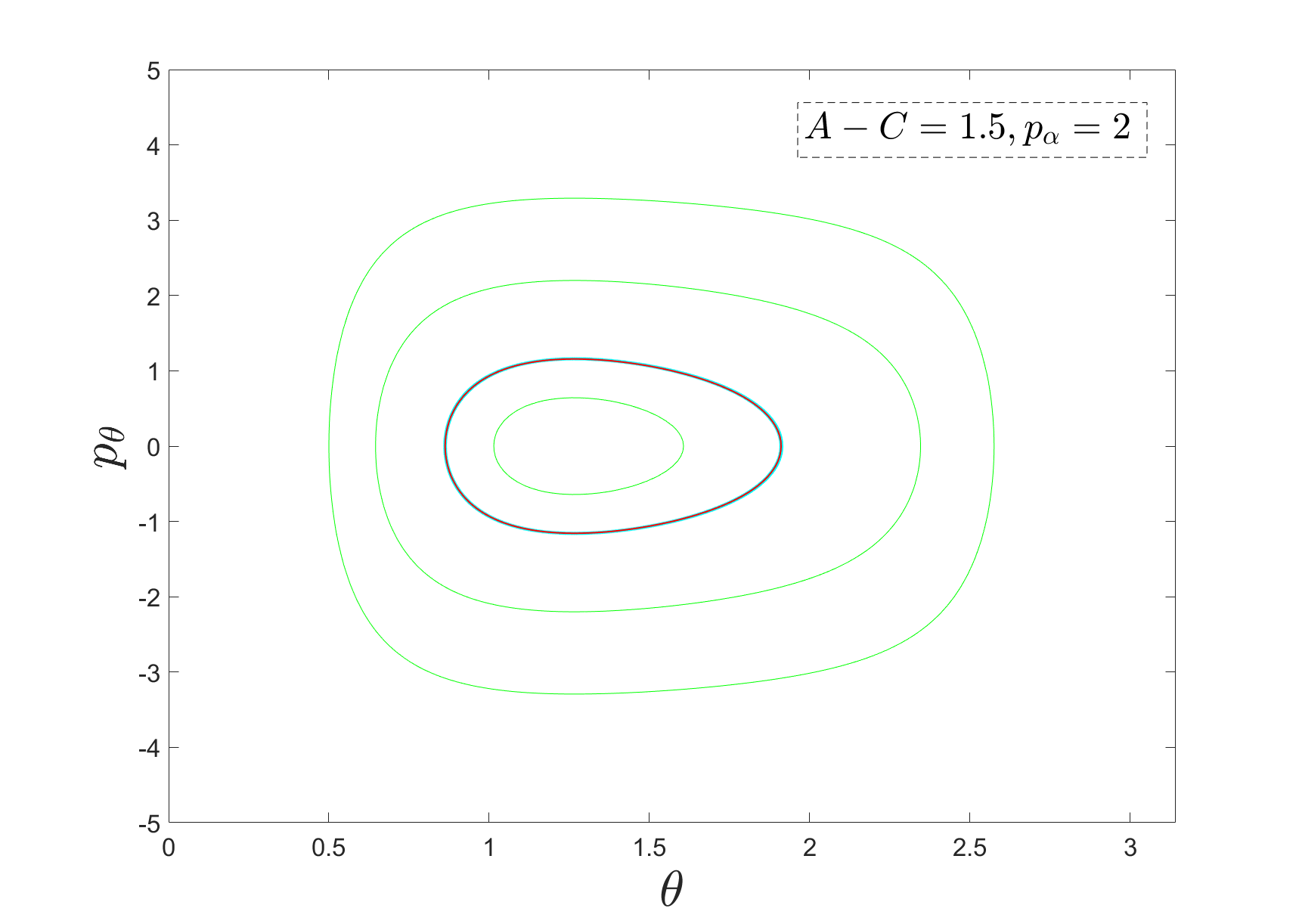}\label{p_phi_4}} 
\caption{Trajectories of the exact (cyan) and the averaged  (red) systems.}
\label{projection_P}
\end{figure}

We numerically calculate the trajectories of the exact Hamiltonian system $H$ (\ref{original_H}) with initial conditions $p_{\theta}(0)=1$, $\theta(0)=\pi/2$ and $\alpha(0)=0.5$  for all examples respectively.  Projections onto the $\theta-p_{\theta}$ plane of these trajectories are shown in cyan in Fig. \ref {projection_P}. The trajectories of the averaged system with the same initial conditions are shown in red in the same figure. One can see that the corresponding trajectories agree very well. The time-step of the plot is $0.001$.
%The behaviour of  $p_{\alpha}$ and of the averaged Hamiltonian $\bar H$ along trajectories of the exact system with examples of vibrations 1,2 and 3, 4 are shown in Fig.\! \ref{graph_p_alpha} and Fig.\! \ref{graph_bar_H} respectively. 

\begin{figure}[htbp!]
  \centering
\subfloat[]{\includegraphics[width=0.45\textwidth]{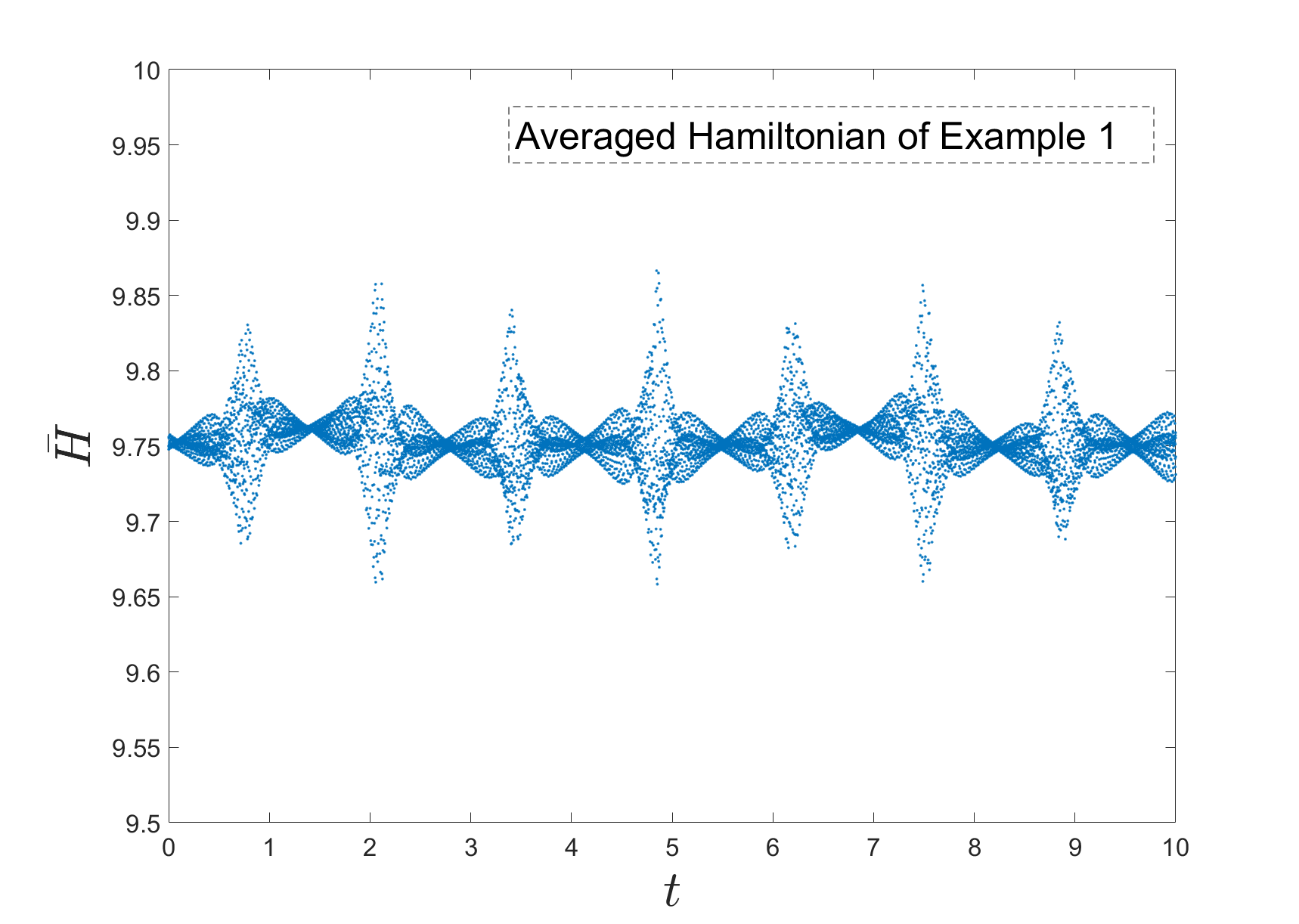}\label{H_1}} 
\subfloat[]{\includegraphics[width=0.45\textwidth]{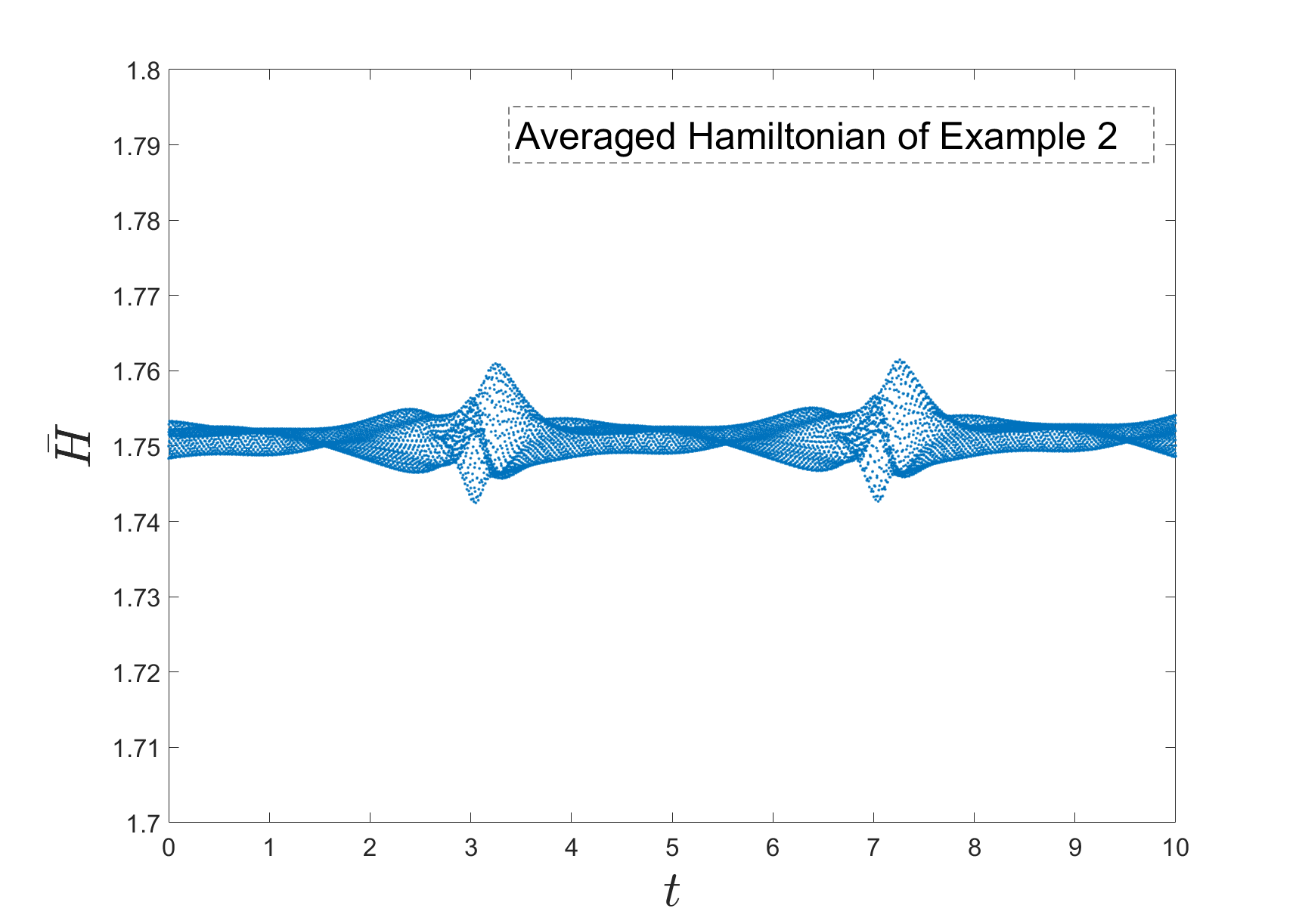}\label{H_2}} 
\\
\subfloat[]{\includegraphics[width=0.45\textwidth]{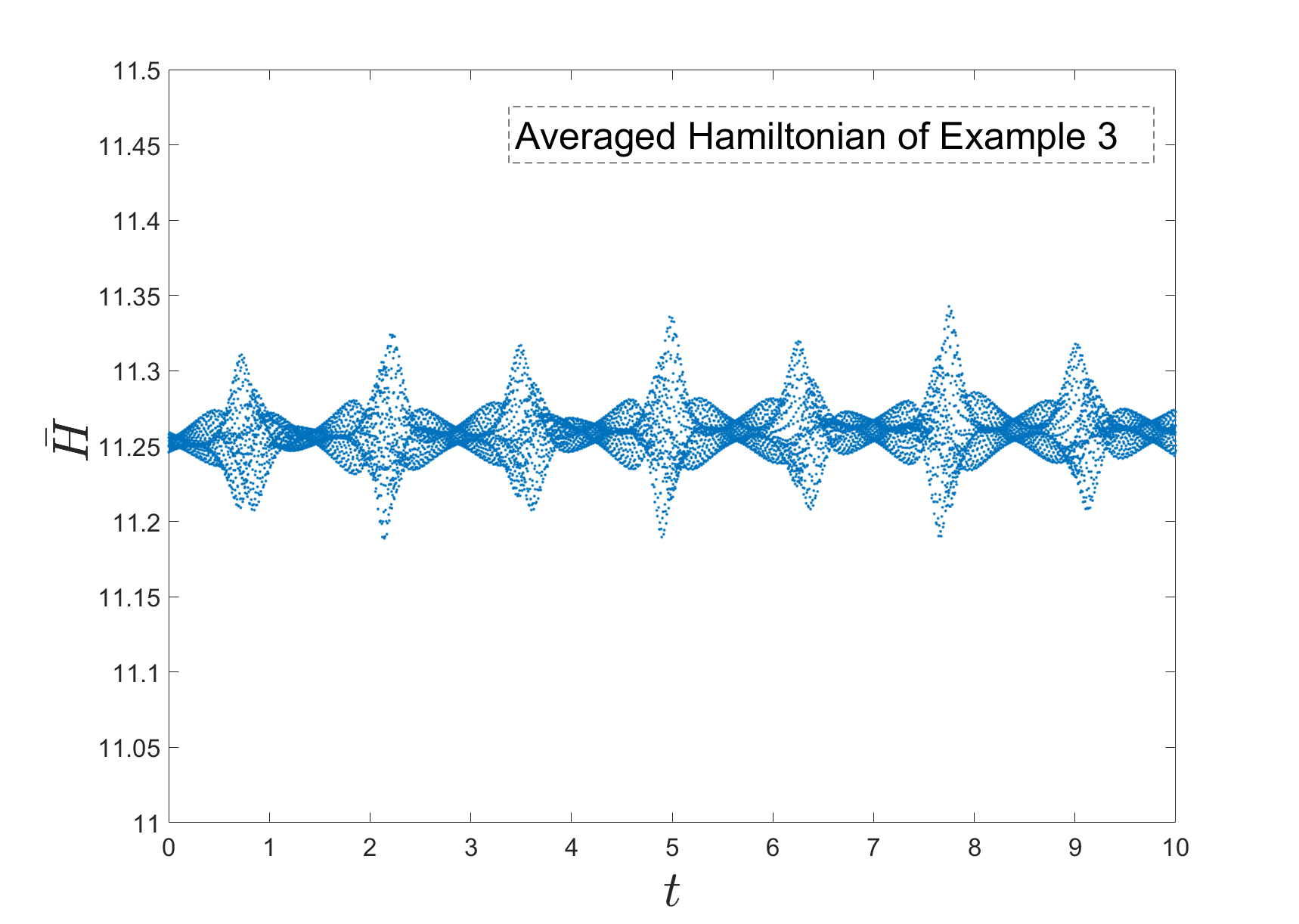}\label{H_3}} 
\subfloat[]{\includegraphics[width=0.45\textwidth]{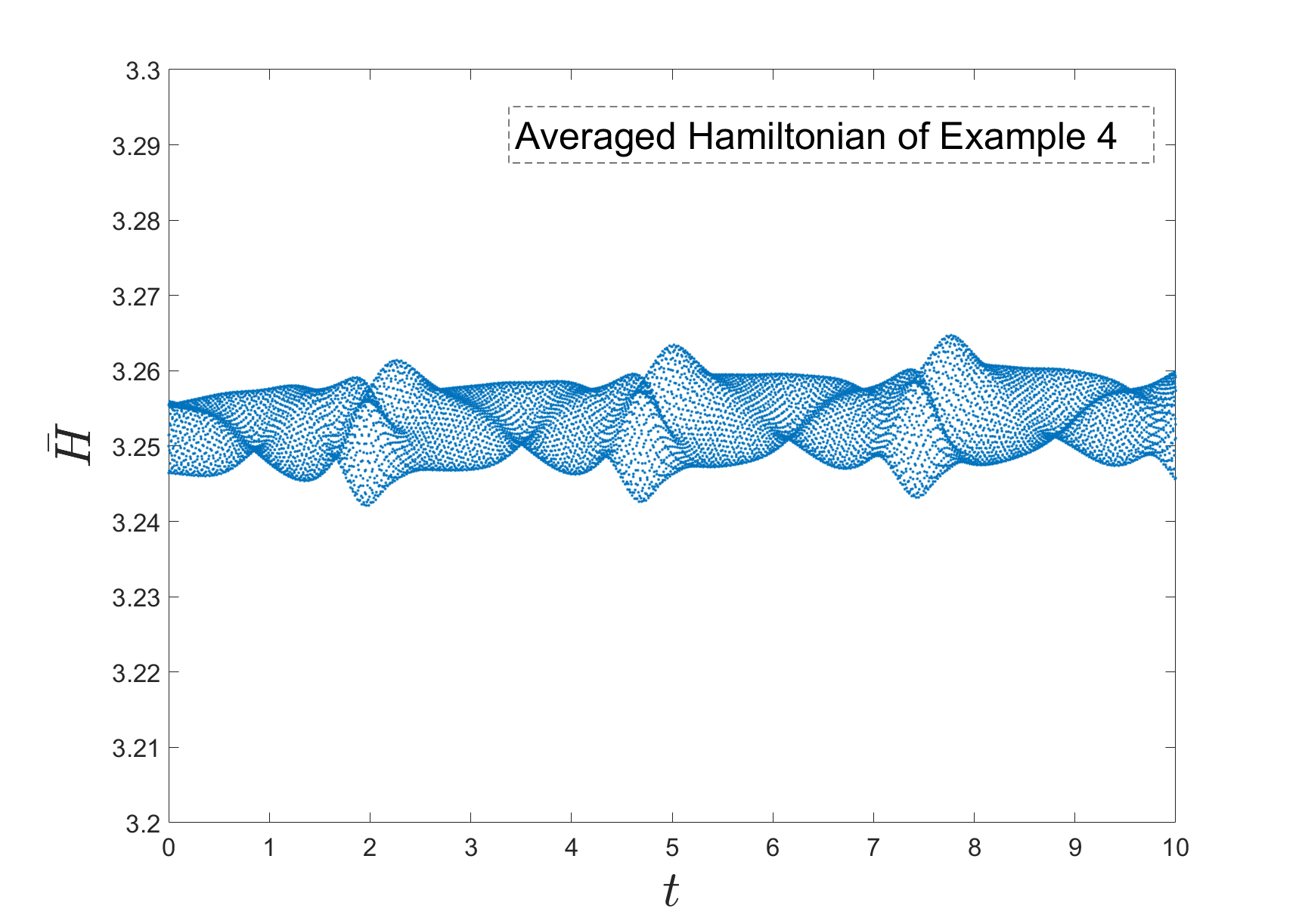}\label{H_4}} 
\caption{Value of $\bar{H}$ vs time.}
\label{graph_bar_H}
\end{figure}

One can see that values of $\bar H$ along the trajectory of the exact system with vibrations in Fig.\ref {graph_bar_H} are subject to only small oscillations, which coincide with our results. 
%Projection onto the plane $\theta, p_{\theta}$ of trajectory of  the exact  system with vibrations  (\ref{example_Q})   and  trajectory of the averaged system with the same initial conditions are shown in Fig.\!  \ref {projection_Q} in blue and red colours respectively.

\section{Conclusion}
We considered a spherical pendulum whose suspension point performs high-frequency spatial vibrations. This system has 2.5 degrees of freedom with canonical phase variables $\theta,\alpha, p_{\theta}, p_{\alpha}$ described in Section \ref{Sec_Hamiltonian}. After averaging over phases of fast vibrations we got a system with 2 degrees of freedom.  Conditions  (\ref{symmetry}) imply that the averaged system has a rotationally symmetry with the cyclic coordinate $\alpha$. Value $p_{\alpha}$ is a first integral of the averaged system. Dynamics of $\theta, p_{\theta}$ is described by a Hamiltonian system with one degree of freedom and can be studied in a phase plane. We described bifurcations of phase portraits of this system. The function $\bar{H}$ (\ref{av_ham}) is its Hamiltonian. KAM theory is applied to demonstrate the relationship between the averaged and the exact system. The dynamical characteristics of the exact problem are similar to the dynamical characteristics in the averaged system.

%\section*{Acknowledgements}

%The authors express their gratitude to Prof. Anatoly Neishtadt for suggestions on the topic and discussions. Kaicheng Sheng thanks the National Natural Science Foundation of China (NSFC) for the support of this research (Grant:  12371192 \& 12271300). 

\bigskip

\section*{Declarations}

\subsection*{Ethics approval and consent to participate}

Not applicable.

\subsection*{Consent for publication}

Not applicable.

\subsection*{Availability of data and materials}

This manuscript does not report data generation or analysis. All of the material is owned by the authors and no permissions are required.

\subsection*{Competing interests}

I declare that the authors have no competing interests as defined by Springer, or other interests that might be perceived to influence the results or discussion reported in this paper.

\subsection*{Funding}

National Natural Science Foundation of China (NSFC), 12371192

National Natural Science Foundation of China (NSFC), 12271300

\subsection*{Authors' contributions}

The authors write the paper together. Sheng performed the analytical work and Luo performed the numerical work. 

\subsection*{Acknowledgements}

The authors express their gratitude to Prof. Anatoly Neishtadt for suggestions on the topic and discussions.

%\medskip
%{ 
%{\bf Data availability}

%\medskip

%Data sharing does not apply to this article as no new data were created or analyzed in this study.

%}
%\medskip
 
%{\bf Conflict of interest statement}

%\medskip
%The authors declare that there are no conflicts of interest, we do not have any possible conflicts of interest.

\normalem
\bibliographystyle{plainnat}  
\bibliography{reference} 

\end{document}